\documentclass{amsart}

\newtheorem{theorem}{Theorem}
\newtheorem{proposition}[theorem]{Proposition}
\newtheorem{lemma}[theorem]{Lemma}
\newtheorem{corollary}[theorem]{Corollary}

\numberwithin{equation}{section}

\newcommand{\Ob}[1]{\mathcal O({1-\beta})^{#1}}
\newcommand\T{\mathcal T}
\newcommand\omegabar{\overline\omega}

\begin{document}

\title[A quadratic approximation to the Sendov radius]{A quadratic 
approximation to the Sendov radius near the unit circle}
\author{Michael J. Miller}
\address{Department of Mathematics\\Le Moyne College\\Syracuse, New 
York 13214}
\email{millermj@mail.lemoyne.edu}
\subjclass{Primary 30C15}
\keywords{Sendov, Ilieff, Ilyeff}
\thanks{30-Sep-2003} 

\begin{abstract}
 Define $S(n,\beta)$ to be the set of complex polynomials of
degree $n\ge2$ with all roots in the unit disk and at least one root 
at $\beta$.  For a polynomial~$P$, define $|P|_\beta$ to be the 
distance between $\beta$ and the closest root of the derivative~$P'$.
Finally, define $r_n(\beta)=\sup \{ |P|_\beta : P \in S(n,\beta) \}$. 
In this notation, a conjecture of Bl. Sendov claims that
$r_n(\beta)\le1$.

In this paper we investigate Sendov's conjecture near the unit circle,
by computing constants $C_1$ and $C_2$ (depending only on $n$) such
that $r_n(\beta)\sim1+C_1(1-|\beta|)+C_2(1-|\beta|)^2$ for $|\beta|$
near $1$.  We also consider some consequences of this approximation.
\end{abstract}

\maketitle

\section{Introduction}  

In 1962, Sendov conjectured that if a polynomial (with complex 
coefficients) has all its roots in the unit disk, then within one unit of
each of its roots lies a root of its derivative.  More than 50 papers 
have been published on this conjecture,   but it has been verified in 
general only for  polynomials of degree at most 8 \cite{Bro}.

Let $n\ge 2$ be an integer and let $\beta$ be a complex number of
modulus at most~$1$.  Define $S(n,\beta)$ to be the set of polynomials 
of degree $n$ with complex coefficients, all roots in the unit disk 
and at least one root at $\beta$.  For a polynomial~$P$, define
$|P|_\beta$ to be the  distance between $\beta$ and the closest root 
of the derivative~$P'$. Finally, define $r_n(\beta)=\sup \{ |P|_\beta 
: P \in S(n,\beta) \}$, and note that $r_n(\beta)\le2$ (since by the
Gauss-Lucas Theorem \cite[Theorem 6.1]{Mar} all roots of each $P'$
are also in the unit disk, and so each $|P|_\beta\le2$). In this
notation, Sendov's conjecture claims simply that $r_n(\beta)\le1$.

In estimating $r_n(\beta)$, we will assume without loss of
generality (by rotation) that $0\le\beta\le1$.   It is already known
that $r_2(\beta)=(1+\beta)/2$ and that 
\begin{equation*}
   r_3(\beta)=[3\beta+(12-3\beta^2)^{1/2}]/6
\end{equation*}
\cite[Theorem 2]{Rah}, that $r_n(0)=(1/n)^{1/(n-1)}$ \cite[Lemma~4 
and $p(z)=z^n-z$]{BRS-2}, that $r_n(1)=1$ \cite[Theorem 1]{Rub}, 
and that $r_n(\beta)\le\min(1.08332, 1+0.72054/n)$ \cite[Corollary 1 
and equations (3)]{BRS-1}.

Since $r_n(1)=1$, an obvious place to look for counterexamples to 
Sendov's conjecture is in a neighborhood of $\beta=1$.  This has 
already been done in \cite[Theorem~3]{Mil-2} and \cite{VZ}, where a 
linear upper bound on $r_n(\beta)$ suffices to verify the Sendov
conjecture  if $\beta$ is sufficiently close to $1$.  Unfortunately,
having only an  upper bound leaves many interesting questions about 
the conjecture unanswered.  In this paper we investigate Sendov's
conjecture much more thoroughly near $\beta=1$, by providing a
quadratic approximation to $r_n(\beta)$ with

\begin{theorem}\label{1}
Let $n\ge3$,  let $k$ be the largest integer such that $k \le (n+1)/3$
and let 

\begin{align*}
u_1 &= \cos{\dfrac{2\pi k}{n+1}},\qquad
                       u_2 = \cos{\dfrac{2\pi (k+1)}{n+1}},\\
D_1 &= {\dfrac{-2 u_1 u_2-1}{2(1-u_1)(1-u_2)}},\qquad
                       D_2 = {\dfrac{-1}{2(1-u_1)(1-u_2)}},\\
D_3 &= (-1-4D_1-3D_1^2+2D_2^2)/2,\\
D_4 &= (3D_1-4D_2+3D_1^2-2D_1D_2-6D_2^2)/2,\\
D_5 &= (2+4D_1+5D_2+2D_1^2+4D_1D_2+3D_2^2)/2, \\
D_6 &= (2D_2+2D_1D_2+3D_2^2)/2 \quad{\rm and}\\
D   &= D_3 n+D_4+D_5/n+D_6/n^2.
\end{align*}

If $n=3$ or $n=5$ then let $\alpha=3/2$; otherwise let $\alpha=2$.
If $n=5$ then let $\Delta=7/225$; otherwise let $\Delta=0$. Then
for $\beta$ sufficiently close to $1$, we have
\begin{equation*}
   r_{n+1}(\beta) = 1+(D_1+D_2/n)(1-\beta) +(D+\Delta)(1-\beta)^2+
    \Ob{\alpha+1}. 
\end{equation*}
\end{theorem}

Before proving this theorem, we will examine some of its
consequences. Our first consequence improves on estimates in \cite{Mil-2}
and \cite{VZ} (by providing a value for the coefficient of the linear
term) with
\begin{corollary}\label{2} 
For all $n\ge2$ we have $r_n(\beta)\le 1-(3/10)(1-\beta)+\Ob{2}.$
\end{corollary}
\begin{proof}
Recall that for $2\le n\le3$ we have formulas for $r_n(\beta)$, and so
the result for those values of $n$ follows from the Taylor series of
these formulas at $\beta=1$.  As we will show in part 6 of Lemma
\ref{8}, the quantity $D_1+D_2/n\le-3/10$ for all $n\ge 3$, and so the
rest of Corollary \ref{2} follows from
Theorem~\ref{1}. 
\end{proof}

As we will show in part 6 of Lemma \ref{8}, at $n=4$ we have 
$D_1+D_2/n=-3/10$, so Corollary \ref{2} provides the smallest possible 
linear upper bound for $r_n(\beta)$ that is independent of~$n$.

A second consequence of Theorem \ref{1} shows that the result of
\cite[Theorem 3]{Mil-2} is the best possible (in the sense that $1/3$
cannot be replaced by a larger number), with
\begin{corollary}\label{3}
There exist constants $K_n>0$ with $\lim_{n \to \infty}K_n=1/3$ such 
that 
\begin{equation*}
   r_{n+1}(\beta)= 1-K_n(1-\beta)+\Ob{2}.
\end{equation*}
\end{corollary}
\begin{proof}
Choose $K_n=-(D_1+D_2/n)$ and note that by Theorem \ref{1} we have 
$r_{n+1}(\beta)=1-K_n(1-\beta)+\Ob{2}$.  As we shall see in parts 
5 and 6 of Lemma \ref{8}, the quantity $D_1+D_2/n$ is negative and tends to 
$-1/3$. 
\end{proof}

Recall that $r_n(0)=(1/n)^{1/(n-1)}$.  This quantity is increasing in
$n$, so it is tempting to conjecture that for all fixed $\beta$ the
quantity $r_n(\beta)$ is increasing in $n$.  Indeed, the graphs in
\cite[figure 4.8]{Mil-1} provide some evidence of this for $n=4$, $6$,
$8$, $10$, and $12$.  Unfortunately, this conjecture is false, as is
shown by
\begin{corollary}\label{4}
For $\beta$ sufficiently close to $1$ we have $r_6(\beta)<r_4(\beta)$.
\end{corollary}
\begin{proof}
By Theorem \ref{1} and the constants we will compute at the beginning of 
section~2, we know that $r_4(\beta)=1-(1/3)(1-\beta)+\Ob{2}$ and
that $r_6(\beta)=1-(11/30)(1-\beta)+\Ob{2}$, and the conclusion 
follows.
\end{proof} 

We will verify Theorem \ref{1} by proving the following three propositions:

\begin{proposition}\label{5}
Assume the notation of Theorem \ref{1}. Then for all polynomials $P \in
S(n+1, \beta)$, we have 
\begin{equation*} 
   |P|_\beta\le 1+(D_1+D_2/n)(1-\beta)+ (D+\Delta)(1-\beta)^2+
   \Ob{\alpha+1}.
\end{equation*}
\end{proposition}

\begin{proposition}\label{6}
There are polynomials $P \in S(6,\beta)$ with 
\begin{equation*}
   |P|_\beta= 1-(11/30)(1-\beta)+(29/450)(1-\beta)^2+\Ob{5/2}.
\end{equation*}
\end{proposition}

\begin{proposition}\label{7}
Assume the notation of Theorem \ref{1}.  Then there are real
polynomials $P \in S(n+1,\beta)$ with
\begin{equation*}
   |P|_\beta= 1+(D_1+D_2/n)(1-\beta) +D(1-\beta)^2+\Ob{\alpha+1}.
\end{equation*}
\end{proposition}

From the definition of $D$ in Theorem \ref{1} and the constants we
will compute at the beginning of section 2, for $n=5$ we have
$D_1+D_2/n=-11/30$ and $D+\Delta=29/450$, so Propositions \ref{5} and
\ref{6} together imply that Theorem \ref{1} is true for $n=5$.  Note
that for $n\ne5$ we have $\Delta=0$, so Propositions \ref{5} and
\ref{7} taken together imply that Theorem \ref{1} is true for $n \ne
5$.

In \cite{Mil-3} it was proved that if $n=5$ and if $\beta$ is
sufficiently close to~$1$ then maximal polynomials in $S(n+1,\beta)$
(those for which $|P|_\beta=r_{n+1}(\beta)$) must be nonreal.  Taken
together, Theorem \ref{1} and Proposition \ref{7} provide strong
evidence that this is true only for $n=5$ (although it is conceivable
that this could fail for higher-order approximations).


\section{Preliminaries}

We begin by computing some values (that we will subsequently need) for
the constants that appear in Theorem~\ref{1}, obtaining:
\begin{center}
$\begin{array}{cccccc}
n & u_1 & u_2 & D_1 & D_2 & D_1+D_2/n \\
3 & 0 & -1 & -1/4 & -1/4 & -1/3 \\
4 & {\dfrac{-1+\sqrt5}{4}} & {\dfrac{-1-\sqrt5}{4}} & -1/5 & -2/5 & -3/10 \\
5 & -1/2 & -1 & -1/3 & -1/6 & -11/30 \\
6 & -0.2225 & -0.9010 & -0.3014 \\
7 & 0 & -0.7071 & -0.2929 & -0.2929 & -0.3347 \\
9 & -0.3090 & -0.8090 \\
10 & -0.1423 & -0.6549 & -0.3138
\end{array}$
\end{center}

We next establish some relationships between these constants with
\begin{lemma}\label{8} Assume the notation of Theorem \ref{1}. Then 
\begin{enumerate}
   \item $u_2< -1/2\le u_1$, and $u_1\le 0$ for $n \ne 4$, and 
$u_2>-1$ for $n\ne3,5$,
   \item $u_1+u_2<0$ and $u_1u_2>-1$,
   \item $2n u_1+n+1\ge1$ and $2n u_2+n+1<0$,
   \item $D_1<0$ and $D_2<0$,
   \item $\lim_{n \to\infty} D_1+D_2/n=-1/3$,
   \item $-1<D_1+D_2/n\le-3/10$, with equality only at $n=4$, and
   \item $1+(1+D_1+D_2)(u_i-1)-D_2(2u_i^2-2)=0$ for $i=1$ and $i=2$.
\end{enumerate}
\end{lemma}
\begin{proof}
From the definition of $k$ in Theorem \ref{1}, the relationship
between $k$ and $n$ depends on the residue of $n$ modulo 3.  For
increasing values of $n$ in each of the three residue classes, the
sequence $k/(n+1)$ increases to (or is equal to) $1/3$ and the
sequence $(k+1)/(n+1)$ strictly decreases to $1/3$, so the values of
$u_1$ decrease to (or are equal to) $-1/2$ and the values of $u_2$
strictly increase to $-1/2$.  Since the values of $u_1$ decrease (or
remain constant) in each residue class, and since $u_1\le0$ for $n=3$,
5 and 7 then $u_1\le0$ for all $n\ne4$.  Since the values of $u_2$
strictly increase in each residue class, and since $u_2>-1$ for $n=4$
and $u_2=-1$ for $n=3$ and $n=5$, then $u_2>-1$ for $n\ne3,5$.  This
completes the proof of part 1 of the lemma.

For $n=4$, we have $u_1+u_2=-1/2$ and $u_1u_2=-1/4$.  For $n\ne4$ we
have from part 1 that $u_2<u_1\le0$, and part 2 of the lemma follows 
trivially.

Since $u_1\ge-1/2$, then $2nu_1+n+1\ge1$.  For $n=3$, $4$ and $5$ we
have $\hbox{(k+1)/(n+1)}\le1/2$.  Since in each residue class this 
quotient strictly decreases to $1/3$ then for all $n\ge3$ we have 
$2\pi(k+1)/(n+1)\in (2\pi/3,\pi]$.  Now $\cos x \le 
1/2-3x/(2\pi)$ on this interval, and from the definition of $k$ in 
Theorem \ref{1} we know that $k\ge(n-1)/3$, so 
\begin{equation*}
   u_2=\cos\frac{2\pi(k+1)}{n+1}\le \frac{1}{2}-\frac{3(k+1)}{n+1} \le 
   \frac{1}{2}-\frac{n+2}{n+1}< -\frac{n+1}{2n}
\end{equation*}
which completes the proof of part 3 of the lemma.

At $n=4$, we have $D_1=-1/5$ and $D_2=-2/5$.  For $n\ne4$ we know 
from part 1 of Lemma \ref{8} that $u_2<u_1\le0$ so from the definitions of
$D_1$ and $D_2$ in Theorem \ref{1} we see that $D_1<0$ and $D_2<0$.  This
completes the proof of part 4 of the lemma.

As $n$ tends to infinity, $u_1$ and $u_2$ tend to $-1/2$, so $D_1$ tends 
to $-1/3$ and $D_2$ is bounded.  This completes the proof of part 5 
of the lemma.

By part 2 of Lemma \ref{8} we have $u_1+u_2<0$ and $u_1u_2>-1$.  Since by 
part 4 of Lemma \ref{8} we know that $D_2<0$ then 
\begin{equation*}
   D_1+D_2/n > D_1+D_2=-\frac{1+u_1u_2}{1+u_1u_2-(u_1+u_2)} > -1.
\end{equation*}
From part 1 of Lemma \ref{8} we know that $u_2<-1/2\le u_1$, so by computing 
the partial derivatives of $D_1$ we see that $\partial D_1/\partial
u_1>0$  and $\partial D_1/\partial u_2\le0$.  Since in each residue
class $u_1$ decreases to $-1/2$ and $u_2$ increases to $-1/2$, then in
each residue class $D_1$ decreases to~$-1/3$.  At  $n=5$, $6$ and $10$
we have $D_1<-3/10$, and hence $D_1+D_2/n<D_1<-3/10$ for all
$n\ge3$ except possibly $n=3$, $4$ and $7$. Checking the values of 
$D_1+D_2/n$ (computed at the beginning of section 2) for these 
exceptional values completes the proof of part 6 of the lemma.

Expressing $D_1$ and $D_2$ in terms of $u_1$ and $u_2$ and simplifying
the result verifies part 7, and thus completes the proof of Lemma \ref{8}.
\end{proof}

We now estimate the size of the coefficients of $P'$ with
\begin{proposition}\label{9}
Suppose that $P \in S(n+1, \beta)$ with $P'$ monic and $|P|_\beta
\ge \beta$.  Let 
$P'(z)=\prod_{j=1}^n (z-\zeta_j)=z^n+a_{n-1}z^{n-1}+\dots+a_0$. Then 
\begin{enumerate}
    \item each $\Re[\zeta_j]=\Ob{}$ and each $\Im[\zeta_j]=\Ob{1/2}$, 
    \item each $a_{n-k}=\Ob{k/2}$,
    \item for $k$ odd, each $\Re[a_{n-k}]=\Ob{(k+1)/2}$, and
    \item for $k$ even, each $\Im[a_{n-k}]=\Ob{(k+1)/2}$.
\end{enumerate}
\end{proposition}
\begin{proof}
Parts 1--3 were proved in \cite[Proposition 4]{Mil-3}.  Part 4 is 
proved similarly to part~3, by noting that each term of $\Im[a_{n-k}]$
is a product of $k$ of the $\Re[\zeta_j]$'s and $\Im[\zeta_j]$'s, and 
that for $k$ even, each term has at least one $\Re[\zeta_j]$, so 
from part 1 of Proposition \ref{9} we have that $\Im[a_{n-k}]=\Ob{(k+1)/2}$. 
\end{proof}

To have $P\in S(n+1,\beta)$ requires that the moduli of the roots 
of $P$ are all at most $1$.  We estimate these moduli with

\begin{proposition}\label{10}
Assume the notation of Theorem \ref{1}.  Let $P$ be a polynomial with 
$P'(z)=z^n+a_{n-1}z^{n-1}+\dots+a_0$ and $P(\beta)=0$.  
Let $z\ne\beta$ be a root of $P$, let $\omega$ be the $(n+1)$th root of $1$ 
that is closest to $z$ and let 
$R=(1-\beta)+a_{n-1}(\omega^n-1)/n+\dots+a_0(\omega-1)$. 
\begin{enumerate}
\item For $0<r\le1$, if each $a_k=\Ob{r}$ then
$|z|^2=1-2\Re[R]+\Ob{2r}$.
\item Suppose that
\begin{align*}
a_{n-1}&=n(1+D_1+D_2)(1-\beta)+\Ob{\alpha},\\
a_{n-2}&=-(n-1)D_2(1-\beta)+\Ob{\alpha}, and\\
a_{n-k}&=\Ob{\alpha} \text{\quad for $k\ge3$}\\
\end{align*}
and define
\begin{align*}
\Gamma_2 &= 2(1+D_1+D_2)(D_1-2D_2+nD_2)\text{ and}\\
\Gamma_1 &= -\Gamma_2+(-2-4D_1)n+(1+4D_1-4D_2).
\end{align*}
If $\Re[\omega]=u_i$ for $i=1$ or $i=2$ then
\begin{equation*}
   |z|^{2n+2}=1-2(n+1)\Re[R]+(n+1)(\Gamma_1+\Gamma_2u_i)(1-\beta)^2+ 
   \Ob{\alpha+1}.
\end{equation*}
\end{enumerate}
\end{proposition}

\begin{proof}
Since $\beta=1-(1-\beta)$ then by the binomial theorem
$\beta^k=1-k(1-\beta)+\Ob{2}$.  Since $z$ is a root of $P$ we have
\begin{equation*}
   0=P(z)=\int_\beta^z P'(t)\,dt = \frac{z^{n+1}-\beta^{n+1}}{n+1} + 
   a_{n-1} \frac{z^n-\beta^n}{n} +\dots +a_0(z-\beta),
\end{equation*}
and solving for $z^{n+1}$ gives us 
\begin{equation}\label{e-2.1}
   z^{n+1}=\beta^{n+1}-(n+1)\left[a_{n-1} \frac{z^n-\beta^n}{n}
   +\dots +a_0(z-\beta)\right].
\end{equation}
By hypothesis, as $\beta$ goes to $1$ the $a_k$ all tend to 0 so the 
roots of $P$ tend to the roots of $z^{n+1}-1$, and so the $\omega$ 
appearing in the hypotheses is well-defined.

Now each $\beta^k=1+\Ob{}$, and by the hypothesis of part 1 each 
$a_k=\Ob{r}$.  Putting these estimates into equation \eqref{e-2.1}, we 
see that $z^{n+1}=1+\Ob{r}$.  Then $z=\omega+\Ob{r}$ and so  
$(z^k-\beta^k)/k=(\omega^k-1)/k+\Ob{r}$.  Now note that each 
$a_{n-k}=\Ob{r}$ and that each $\beta^k=1-k(1-\beta)+\Ob{2}$. 
Substituting these estimates into equation~\eqref{e-2.1} gives 
\begin{equation*}
\begin{split}
z^{n+1}&=1-(n+1)(1-\beta)-(n+1)\left[a_{n-1} \frac{\omega^n-1}{n} +\dots 
+a_0(\omega-1) \right]+\Ob{2r} \\
          &=1-(n+1)R+\Ob{2r}. 
\end{split}
\end{equation*}
Note that $R=\Ob{r}$ so
\begin{equation*}
   (1-R)^{n+1}=1-(n+1)R+\Ob{2r}=z^{n+1}+\Ob{2r}
\end{equation*}
so 
$z=\omega(1-R)+\Ob{2r}$ and  hence $|z|^2=z\overline z=
1-2\Re[R]+\Ob{2r}$.  This finishes the proof of  part 1.

From the hypotheses of part 2, we know that $\Re[\omega]=u_i$ for 
$i=1$ or $i=2$. Suppose for the moment that $\Re[\omega]=u_1$ and 
write $\omega=u_1+iv_1$. Since $\omega^{n+1}=1$ then $|\omega|=1$, so 
$\omega^n=\omegabar$ and $\Re[\omega^2]=2u_1^2-1$. Let
$A=\big[-(1+D_1+D_2)+2D_2u_1\big]v_1$.  From part 7 of Lemma \ref{8} we see 
that 
\begin{equation*}
   \Re[1+(1+D_1+D_2)(\omegabar-1)-D_2(\omegabar^2-1)]=0
\end{equation*}
and so using the estimates of the $a_{n-k}$'s given in the hypotheses 
of part 2, we get 
\begin{equation*}
\begin{split}
R &= (1-\beta)+ a_{n-1}\frac{\omegabar-1}{n}+ a_{n-2}\frac{\omegabar^2-1}{ 
n-1}+ \dots + a_0(\omega-1) \\
&=\big[1+(1+D_1+D_2)(\omegabar-1)-D_2(\omegabar^2-1)\big](1-\beta)+
\Ob{\alpha} \\
  &=iA(1-\beta)+\Ob{\alpha}.
\end{split}
\end{equation*}

The hypotheses of part 2 imply that each $a_k=\Ob{}$, so from the 
proof of part 1 with $r=1$ we have $z=\omega(1-R)+ 
\Ob{2}=\omega\big[1-i A(1-\beta)\big]+\Ob{\alpha}$ 
and so
\begin{equation*}
   (z^k-\beta^k)/k=(\omega^k-1)/k+(1-i A\omega^k)(1-\beta)+\Ob{\alpha}.
\end{equation*}

Let $G=n/2-n(1+D_1+D_2)(1-i A\omegabar)+(n-1)D_2(1-i A\omegabar^2)$.
Then from equation \eqref{e-2.1} and the estimates of the $a_k$'s 
given in the hypotheses of part 2 we get 
\begin{equation*}
\begin{split}
z^{n+1}&=1-(n+1)(1-\beta)+\frac{(n+1)n}{2}(1-\beta)^2 \\
         &\qquad-(n+1)\bigg[a_{n-1}\left(\frac{\omega^n-1}{ 
n} +(1-i A\omega^n)(1-\beta)\right) \\
         &\qquad\qquad\qquad\qquad+a_{n-2}\left(\frac{\omega^{n-1}-1}{ 
n-1} +(1-i A\omega^{n-1})(1-\beta)\right) \\
         &\qquad\qquad\qquad\qquad+a_{n-3}\frac{\omega^{n-2}-1}{n-2}+ 
\dots+ a_0(\omega-1)\bigg]+\Ob{\alpha+1} \\
         &=1-(n+1)R+(n+1)G(1-\beta)^2+\Ob{\alpha+1}. 
\end{split}
\end{equation*}
Then since $R=iA(1-\beta)+\Ob{\alpha}$ we have
\begin{equation*}
\begin{split}
|z|^{2n+2}&=z^{n+1}\overline z^{n+1} \\
            &=1-2(n+1)\Re[R]+
(n+1)\big[2\Re[G]+(n+1)A^2\big](1-\beta)^2+\Ob{\alpha+1}. 
\end{split}
\end{equation*} 

Thus to complete the proof of part 2 of Proposition \ref{10} for the
case $\Re[\omega]=u_1$ we need only verify that
$2\Re[G]+(n+1)A^2=\Gamma_1+\Gamma_2u_1$.

Let $D_0=1+D_1+D_2$, so from the definition of $A$ we see that
$A=\break(-D_0+2D_2u_1)v_1$. Note that
$\Re[i\,\omegabar]=\Im[\omega]$.  Then from the definition of $G$ we
have
\begin{equation*}
\begin{split}
\Re[G] &= n/2-nD_0(1-Av_1)+(n-1)D_2(1-2Au_1v_1) \\
    &= n/2-nD_0+(n-1)D_2
            -A\big[n(-D_0v_1+2D_2u_1v_1)-2D_2u_1v_1\big] \\ 
   &= (-n/2-nD_1-D_2)-nA^2+2AD_2u_1v_1
\end{split}
\end{equation*}
so 
\begin{equation}\label{e-2.2}
   2\Re[G]+(n+1)A^2=(-n-2nD_1-2D_2)+(-n+1)A^2+4AD_2u_1v_1.
\end{equation}

Now
\begin{equation*}
   2D_2u_1^2=\frac{-u_1^2}{(1-u_1)(1-u_2)}=D_0u_1+(D_2-D_1)
\end{equation*}
so
\begin{equation*}
\begin{split}
Av_1 &= (-D_0+2D_2u_1)(1-u_1^2) \\
     &= -D_0+2D_2u_1-u_1(-D_0u_1+2D_2u_1^2) \\
     &= -D_0+(D_1+D_2)u_1.
\end{split}
\end{equation*}

Using these two equalities, we see that
\begin{equation*}
\begin{split}
A^2 &= (-D_0+2D_2u_1) \big[-D_0+(D_1+D_2)u_1\big] \\
     &= D_0^2+(-D_0D_1-3D_0D_2)u_1+(D_1+D_2)(2D_2u_1^2) \\
     &= D_0^2-D_1^2+D_2^2-2D_0D_2u_1
\end{split}
\end{equation*}
and
\begin{equation*}
\begin{split}
2AD_2u_1v_1 &= 2D_2u_1[-D_0+(D_1+D_2)u_1] \\
                 &=-2D_0D_2u_1+(D_1+D_2)[D_0u_1+(D_2-D_1)] \\
                &= D_0(D_1-D_2)u_1+(D_2^2-D_1^2).
\end{split}
\end{equation*}
Thus from equation \eqref{e-2.2} we have
\begin{equation*}
\begin{split}
2\Re[G]+(n+1)A^2 &= (-n-2nD_1-2D_2) 
+(-n+1)(D_0^2-D_1^2+D_2^2-2D_0D_2u_1) \\ 
               &\qquad+ 2\big[D_0(D_1-D_2)u_1+(D_2^2-D_1^2)\big] \\
      &= (-1-2D_1-D_0^2+D_1^2-D_2^2)n +(-2D_2+D_0^2-3D_1^2+3D_2^2) \\ 
           &\qquad+ 2D_0u_1(D_1-2D_2+nD_2) \\
     &=\Gamma_1+\Gamma_2u_1.
\end{split}
\end{equation*}

This finishes the proof of part 2 of Proposition \ref{10} for the case 
$\Re[\omega]=u_1$.  Since $D_1$ and $D_2$ are symmetric in $u_1$ and 
$u_2$, swapping $u_1$ and $u_2$ in this proof verifies part 2 of
Proposition \ref{10} for the remaining case $\Re[\omega]=u_2$, and thus
completes the proof of Proposition \ref{10}.
\end{proof}

Finally, consider the linear transformation $\T$ which takes functions
to real numbers via 
\begin{equation}\label{e-2.3}
   \T(f)=\frac{(2nu_1+n+1)f(u_2)-(2nu_2+n+1)f(u_1)}{2(u_1-u_2)}.
\end{equation}
Recall that by Lemma \ref{8} we have $u_1-u_2>0$, $2nu_1+n+1>0$ and 
$2nu_2+n+1<0$, so $\T/n$ is a weighted average.  This implies that 
$\T$ preserves inequalities, in the sense that if $f(u_1)\le
g(u_1)$ and $f(u_2)\le g(u_2)$ then $\T(f)\le\T(g)$.

In the process of analyzing several inequalities, we will need the
following values of the transformation $\T$:
\begin{equation}\label{e-2.4}
\begin{split}
   \T(1) &= n \\
   \T(2+2u)&=n-1 \\
   \T(1+4u+4u^2)&=-[n+2+2(n+1)(u_1+u_2)+4nu_1u_2] \\
                        &=-\frac{n+1+D_1+3nD_1+3D_2}{D_2} \\
   \T\left(\frac{1}{1-u}\right) &= n+nD_1+D_2 \\
   \T\left(\frac{u}{1-u}\right) &= nD_1+D_2 
\end{split}
\end{equation}

We will also use the results of
\begin{lemma}\label{11} For the linear transformation $\T$ defined in 
equation \eqref{e-2.3} we have 
\begin{enumerate}
   \item $\T(1+4u+4u^2)/(n-2)<1/2$ for $n\ne3$, $4$ and $6$, and
   \item $\T(8u^2+8u^3)\ge0$ for all $n$.
\end{enumerate}
\end{lemma}
\begin{proof}
From the formula for $\T(1+4u+4u^2)$ in \eqref{e-2.4} and from 
part 3 of Lemma \ref{8} we have
\begin{align*}
\partial \T(1+4u+4u^2)/\partial u_1 &= -2(2nu_2+n+1)>0 \text{ 
and} \\ 
\partial \T(1+4u+4u^2)/\partial u_2 &= -2(2nu_1+n+1)<0.
\end{align*}
Recall from the proof of Lemma \ref{8} that for each residue class of
$n$ modulo $3$ the values of $u_1$ decrease and the values of $u_2$
increase, so the signs of the partial derivatives above imply that in
each residue class the values of $\T(1+4u+4u^2)$ decrease.  Since
$1+4u+4u^2= (1+2u)^2\ge0$ and since $\T$ preserves inequalities, then
$\T(1+4u+4u^2)\ge0$, so the values of $\T(1+4u+4u^2)/(n-2)$ also
decrease in each residue class.  Using the formula for $\T(1+4u+4u^2)$
in \eqref{e-2.4} and the values of the $u_i$ computed at the beginning
of section 2, we calculate the values of $\T(1+4u+4u^2)/(n-2)$ at
$n=5$, $7$ and $9$, getting respectively $1/3$, $0.4627$ and $0.3372$.
Since they are all less than $1/2$, this proves part~1 of Lemma
\ref{11}.

Since by definition $u_i\ge-1$ then $8u_i^2+8u_i^3=8u_i^2(1+u_i)\ge0$ 
for both $i=1$ and $i=2$, and so part 2 of Lemma \ref{11} follows from our
observation that $\T$ preserves inequalities.
\end{proof}

Finally, we will deal with polynomials that are ``almost'' in $S(n,\beta)$
using

\begin{lemma}\label{12}  
Suppose that $P$ is a polynomial of degree $n$ with all roots in $\{z:
|z|\le1+\Ob{r}\}$, one root at $\beta$, and all other roots bounded
away from $\beta$.  Then there is a polynomial $Q \in S(n,\beta)$ such
that $|Q|_\beta = |P|_\beta +\Ob{r}$.
\end{lemma}
\begin{proof}
If $P \in S(n,\beta)$ then we may take $Q=P$.  If not, then at least
one root of $P$ has modulus greater than $1$.  In this case, let
\begin{equation*}
c=\max\left\{\frac{|z|^2-1}{|z-\beta|^2}: \text{$z$ is a root of $P$ and
$|z|>1$}\right\}
\end{equation*}
Since by hypothesis $|z-\beta|$ is bounded away from $0$ and
$|z|\le1+\Ob{r}$, then $0<c\le\Ob{r}$.  In particular, for $\beta$
sufficiently close to $1$ we have $0<c<1$.

Let $Q$ be the polynomial with roots $\{z-c(z-\beta): \text{$z$ is 
a root of $P$}\}$.  Since the mapping $z \mapsto z-c(z-\beta)$ is a 
contraction of the plane that leaves $\beta$ fixed and moves all roots
of $P$ (and hence $P'$) toward $\beta$ by at most $\Ob{r}$, then 
$Q(\beta)=0$ and $|Q|_\beta=|P|_\beta+\Ob{r}$.  Thus we need only show
that all roots of $Q$ are in the unit disk.

Note that for $t$ real the image of the mapping $t \mapsto
z-t(z-\beta)$ is a line, with $t=0$ mapping to $z$, and $t=1$ mapping 
to $\beta$, and $t=(|z|^2-1)/|z-\beta|^2$ mapping to
\begin{equation*}
z-\frac{|z|^2-1}{|z-\beta|^2}(z-\beta) = 
   z-\frac{z\overline z-1}{\overline z-\beta} = 
   \frac{1-\beta z}{\overline z-\beta}.
\end{equation*}
If $z$ is in the unit disk, then the images of every $t$ between $0$
and $1$ lie on the line between $z$ and $\beta$, hence in the unit
disk.  If $z$ is not in the unit disk, then $|(1-\beta
z)/(z-\beta)|<1$ and so the images of every $t$ between
$(|z|^2-1)/|z-\beta|^2$ and $1$ lie on the line between $(1-\beta
z)/(\overline z-\beta)$ and $\beta$, hence in the unit disk.  Thus for
every root $z$ of $P$, the image of $c$ lies in the unit disk, so all
roots of $Q$ are in the unit disk and so $Q \in S(n,\beta)$. This
completes the proof of Lemma \ref{12}.
\end{proof}


\section{Proof of Proposition \ref{5}}

Take any $P \in S(n+1,\beta)$, assume without loss of generality 
that $P'$ is monic, and write $P'(z)=\prod_{j=1}^n(z-\zeta_j)=  
z^n+a_{n-1}z^{n-1}+\dots+a_0$.  

If $|P|_\beta\le 1+(D_1+D_2/n)(1-\beta)+ (D+\Delta)(1-\beta)^2$, then 
Proposition~\ref{5} is trivially true.  Thus we may assume without loss of
generality that 
\begin{equation}\label{e-3.1}
   |P|_\beta\ge 1+(D_1+D_2/n)(1-\beta)+(D+\Delta)(1-\beta)^2.
\end{equation}
From part 6 of Lemma \ref{8} we have that $D_1+D_2/n>-1$,  and so inequality 
\eqref{e-3.1} implies that $|P|_\beta\ge\beta$ as long as $\beta$ is
sufficiently close to $1$.  Note that $P$ thus satisfies all the 
hypotheses of Proposition \ref{9}.

We begin by estimating some relationships between the coefficients of 
$P'$ with

\begin{lemma}\label{13}
Suppose that $\Im[a_{n-1}]=\Ob{3/2}$ and that each 
\begin{equation*}
   |\zeta_j-\beta|=1+(D_1+D_2/n)(1-\beta)+\Ob{2}.
\end{equation*}
Then
\begin{enumerate}
\item $\Im[a_{n-2}]= (-3/2)\Im[a_{n-3}]+\Ob{5/2}$ and
\item $\Re[a_{n-3}]+2\Re[a_{n-4}]=
(n-2)(1+D_1+D_2/n)(1-\beta)\Re[a_{n-2}]+\Ob{3}$.
\end{enumerate}
\end{lemma}
\begin{proof}
Let each $\zeta_j=x_j+i y_j$ and note that by Proposition \ref{9} we have 
$x_j=\Ob{}$ and $y_j=\Ob{1/2}$.  Note that by hypothesis, $\sum_i y_i= 
-\Im[a_{n-1}]=\Ob{3/2}$ and that each 
\begin{equation*}
   (\beta-x_j)^2+y_j^2 = |\beta-\zeta_j|^2= 1+2(D_1+D_2/n)(1-\beta)+\Ob{2}
\end{equation*}
so solving for $x_j$ gives us
\begin{equation}\label{e-3.2}
   x_j=y_j^2/2-(1+D_1+D_2/n)(1-\beta)+\Ob{2}.
\end{equation}

Note that $\Im[a_{n-3}]=-\sum_{i<j<k} \Im[\zeta_i \zeta_j \zeta_k] 
=\sum_{i<j<k} y_i y_j y_k+\Ob{5/2}$, so
\begin{equation*}
\begin{split}
\Ob{5/2} &= \sum_i y_i \sum_{i<j} y_i y_j= \sum_{i\ne j}y_i^2 y_j
+3\sum_{i<j<k} y_i y_j y_k \\
            &= \sum_{i\ne j} y_i^2y_j+3\Im[a_{n-3}]+\Ob{5/2}
\end{split}
\end{equation*}
and so $\sum_{i\ne j} y_i^2y_j=-3\Im[a_{n-3}]+\Ob{5/2}$. Then using
equation \eqref{e-3.2} we have
\begin{equation*}
\begin{split}
\Im[a_{n-2}] &= \sum_{i<j}\Im[\zeta_i \zeta_j] = \sum_{i \ne j}x_i 
y_j \\
     & = (1/2)\sum_{i\ne j} y_i^2 y_j- (1+D_1+D_2/n)(1-\beta)\sum_{i\ne 
j} y_j+\Ob{5/2} \\
                &= (-3/2)\Im[a_{n-3}]+\Ob{5/2},
\end{split}
\end{equation*}
which completes the proof of part 1 of Lemma \ref{13}.

Let $S$ be the set of triples $(i,j,k)$ of distinct integers from $1$ 
to $n$ with $j<\nobreak k$.   Note that $\Re[a_{n-2}] = \sum_{i<j} 
\Re[\zeta_i\zeta_j]= -\sum_{i<j}y_iy_j+\Ob{2}$ and
$\Re[a_{n-3}]=-\sum_{i<j<k} \Re[\zeta_i\zeta_j\zeta_k]=\sum_S
x_iy_jy_k+\Ob{3}$.  Furthermore, 
\begin{equation*}
   \Ob{3}=\sum_i y_i \sum_{j<k<l} y_jy_ky_l= \sum_Sy_i^2y_jy_k+
   4\sum_{i<j<k<l} y_iy_jy_ky_l,
\end{equation*}
so 
\begin{equation*}
\begin{split}
\Re[a_{n-4}]&=\sum_{i<j<k<l}\Re[\zeta_i\zeta_j\zeta_k\zeta_l]= 
\sum_{i<j<k<l} y_iy_jy_ky_l+\Ob{3} \\
                 &= (-1/4)\sum_Sy_i^2y_jy_k+\Ob{3}.
\end{split}
\end{equation*}
Then using equation \eqref{e-3.2} we have
\begin{equation*}
\begin{split}
\Re[a_{n-3}]+2\Re[a_{n-4}]
     &=\sum_S \big( x_i-y_i^2/2 \big) y_jy_k+\Ob{3} \\
     &=-(1+D_1+D_2/n)(1-\beta)(n-2)\sum_{j<k} y_jy_k +\Ob{3} \\
     &=(n-2)(1+D_1+D_2/n)(1-\beta)\Re[a_{n-2}]+\Ob{3},
\end{split}
\end{equation*}
which completes the proof of Lemma \ref{13}.
\end{proof}

We now establish a lower bound on $\Re[a_{n-4}]$ with
\begin{lemma}\label{14}
Suppose that
\begin{align*}
\Im[a_{n-1}] &= \Ob{\alpha}, \\
\Re[a_{n-2}] &= -(n-1)D_2(1-\beta)+\Ob{\alpha},\text{ and} \\
\Im[a_{n-3}] &= \Ob{\alpha}.
\end{align*}
If $n=5$ then define $\delta=-1/15$; otherwise define $\delta=0$.
Then 
\begin{equation*}
   \Re[a_{n-4}]\ge\delta(1-\beta)^2+\Ob{\alpha+1}.
\end{equation*}
\end{lemma}
\begin{proof}
Let each $\zeta_j=x_j+iy_j$ and recall by Proposition \ref{9} that 
$x_j=\Ob{}$ and  $y_j=\Ob{1/2}$.  Let $F(y)= \prod_{i=1}^n (y+y_i)=
y^n+b_{n-1}y^{n-1}+\dots+b_0$.  Note that 
\begin{equation*}
\begin{split}
\Re[a_{n-4}] &= 
\sum_{i<j<k<l}\Re[\zeta_i\zeta_j\zeta_k\zeta_l]=      
\sum_{i<j<k<l} y_iy_jy_ky_l+\Ob{3} \\
           &=b_{n-4}+\Ob{3}
\end{split}
\end{equation*}
and that by hypothesis 
\begin{align*}
b_{n-1} &= \sum_i y_i = \sum_i\Im[\zeta_i]= -\Im[a_{n-1}] \\
           &= \Ob{\alpha}, \\
b_{n-2} &= \sum_{i<j} y_i y_j = -\sum_{i<j}\Re[\zeta_i \zeta_j]+\Ob{2} 
=-\Re[a_{n-2}]+\Ob{2} \\
           &= (n-1)D_2(1-\beta)+\Ob{\alpha}, \text{ and} \\
b_{n-3} &= \sum_{i<j<k} y_iy_jy_k = -\sum_{i<j<k} \Im[\zeta_i \zeta_j 
\zeta_k]+\Ob{5/2} =\Im[a_{n-3}]+\Ob{5/2} \\ 
        &= \Ob{\alpha}.
\end{align*}
Let
\begin{equation*}
\begin{split}
f(y) &=F^{(n-4)}(y)\\
     &=\frac{n!}{24}y^4+\frac{(n-1)!}{6}b_{n-1}y^3+
        \frac{(n-2)!}{2}b_{n-2}y^2
        +(n-3)!\,b_{n-3}y+ (n-4)!\,b_{n-4}.
\end{split}
\end{equation*}

Now by definition $F$ has all real roots, hence by Rolle's Theorem
(from elementary calculus) so does $f$. Then the ``reverse'' of f
defined by $y^4f(1/y)= (n-4)!b_{n-4}y^4+\dots+n!/24$ has all real
roots, so by Rolle's theorem so does the reverse's second derivative
\begin{equation*}
   12(n-4)!b_{n-4}y^2+6(n-3)!b_{n-3}y+(n-2)!~b_{n-2}.
\end{equation*}
Since this quadratic has all real roots then its discriminant is 
nonnegative, so
\begin{equation*}
   [6(n-3)!b_{n-3}]^2-48(n-2)!(n-4)!b_{n-2}b_{n-4}\ge0.
\end{equation*}
Using our estimates of the $b_{n-k}$'s (including $b_{n-4}=\Ob{2}$), 
this implies that $-D_2(1-\beta)b_{n-4}\ge\Ob{2\alpha}$ and so 
$b_{n-4}\ge\Ob{2\alpha-1}$.  Now for $n\ne3,5$ we have $\alpha=2$ and 
so $\Re[a_{n-4}]=b_{n-4}+\Ob{3}\ge\Ob{3}$, which finishes the 
proof of Lemma \ref{14} for these values of $n$.

Lemma \ref{14} is trivially true for $n=3$, since then
$\Re[a_{n-4}]\equiv0\ge\Ob{5/2}$.

Finally, for $n=5$ we have that
\begin{equation*}
   f(y)=5y^4+4b_{n-1}y^3+3b_{n-2}y^2+2b_{n-3}y+b_{n-4}
\end{equation*}
has all real roots, hence by Rolle's theorem so does its derivative
\begin{equation*}
   f'(y)=20y^3+12b_{n-1}y^2+6b_{n-2}y+2b_{n-3}.
\end{equation*}
A classical result (see e.g. \cite[p.289]{Usp}) states that if a 
cubic polynomial $ax^3+bx^2+cx+d$ has all real roots then its 
discriminant is nonnegative, so
\begin{equation*}
   18abcd-4b^3d+b^2c^2-4ac^3-27a^2d^2\ge0.
\end{equation*}
Applying this to $f'(y)$, we have
\begin{equation*}
   -4[20][6b_{n-2}]^3-27[20]^2[2b_{n-3}]^2\ge\Ob{4},
\end{equation*}
which implies that $2b_{n-2}^3+5b_{n-3}^2\le\Ob{4}$.
Since for $n=5$ we have $D_2~=~-1/6$, then by hypothesis
 $b_{n-2}=(-2/3)(1-\beta)+\Ob{3/2}$, and so
\begin{equation*}
\begin{split}
b_{n-3}^2 &\le (-2/5)b_{n-2}^3+\Ob{4} \\
          &= (16/135)(1-\beta)^3+\Ob{7/2}.
\end{split}
\end{equation*}

We also have that the first derivative of the reverse of $f$
\begin{equation*}
   4b_{n-4}y^3+6b_{n-3}y^2+6b_{n-2}y+4b_{n-1}
\end{equation*}
has all real roots, so applying our classical result gives
\begin{equation*}
   [6b_{n-3}]^2[6b_{n-2}]^2-4[4b_{n-4}][6b_{n-2}]^3\ge\Ob{6}.
\end{equation*}
Dividing this by $144b_{n-2}^2$ and recalling that $b_{n-2}=
(-2/3)(1-\beta)+\Ob{3/2}$ yields
\begin{equation*}
   9b_{n-3}^2+16(1-\beta)b_{n-4}\ge\Ob{7/2}.
\end{equation*}
Combining these two inequalities implies that for $n=5$ we have
\begin{equation*}
\begin{split}
\Re[a_{n-4}] &= b_{n-4}+\Ob{3} \\
        &\ge \frac{-9b_{n-3}^2}{16(1-\beta)}+\Ob{5/2} \\
        &\ge (-1/15)(1-\beta)^2+\Ob{5/2}.
\end{split}
\end{equation*}
This completes the proof of Lemma \ref{14}.
\end{proof}

We now begin the proof of Proposition \ref{5}.  Our first step will be to
show that $|P|_\beta\le 1+(D_1+D_2/n)(1-\beta)+\Ob{2}$.  Recall that 
$P$ satisfies the hypotheses of Proposition \ref{9}, so each
$a_{n-k}=\Ob{k/2}$.  Let $\omega\ne 1$ be any $(n+1)$st root of $1$ and 
let $z$ be the root of $P$ (so $|z|\le1$) closest to $\omega$.  Then 
in Proposition \ref{10} we have 
\begin{equation*}
\begin{split}
R &= (1-\beta)+a_{n-1}(\omega^n-1)/n+\dots+a_0(\omega-1) \\
   &= a_{n-1}(\omega^n-1)/n+\Ob{}
\end{split}
\end{equation*}
and so by part 1 of Proposition \ref{10} with
$r=1/2$, we have 
\begin{equation*}
   |z|^2=1-2\Re[a_{n-1}(\omega^n-1)/n]+\Ob{}.
\end{equation*}

Since $|z|\le1$ and $\omega^n=\overline\omega$, this implies 
that $\Re[a_{n-1}(\overline\omega-1)] \ge \Ob{}$. 
Expanding the product and noting that by Proposition \ref{9} we have
$\Re[a_{n-1}]=\Ob{}$, we get that $\Im[a_{n-1}]\Im[\omega]\ge\Ob{}$. 
Choosing $\omega$ nonreal and repeating this argument with $\omegabar$ 
substituted for $\omega$ provides that 
$\Im[a_{n-1}]\Im[\omegabar]\ge\Ob{}$ and so $\Im[a_{n-1}]=\Ob{}$. Thus
we have $a_{n-1}=\Ob{}$.

Recall that each $a_{n-k}=\Ob{k/2}$, so we know now that each 
$a_{n-k}=\Ob{}$. Since $\omega^{n-k}=\overline\omega^{k+1}$, by part 1
of Proposition \ref{10} with $r=1$ we have
\begin{equation*}
   |z|^2=1-2\Re\left[(1-\beta)+a_{n-1}\frac{\overline\omega-1}{ 
   n}+a_{n-2}\frac{\overline\omega^2-1}{
   n-1}+ a_{n-3}\frac{\overline\omega^3-1}{n-2}\right]+\Ob{2}.
\end{equation*}

Since $|z|\le1$ this implies that
\begin{equation}\label{e-3.3}
   -\Re\left[a_{n-1}\frac{\overline\omega-1}{ 
   n}+a_{n-2}\frac{\overline\omega^2-1}{
   n-1}+ a_{n-3}\frac{\overline\omega^3-1}{n-2}\right]\le(1-\beta)+\Ob{2}.
\end{equation}

Averaging the expressions obtained by substituting $\omega$ and
$\overline\omega$ into inequality \eqref{e-3.3} and noting that by 
Proposition \ref{9} we have $\Re[a_{n-3}]=\Ob{2}$ we get
\begin{equation}\label{e-3.4}
   \Re[a_{n-1}]\Re\left[\frac{1-\omega}{n}\right]+
   \Re[a_{n-2}]\Re\left[\frac{1-\omega^2}{n-1}\right] \le(1-\beta)+\Ob{2}.
\end{equation}
Let $u=\Re[\omega]$.  Note that since $|\omega|=1$, then
$\Re[\omega^2]=2u^2-1$, so dividing inequality~\eqref{e-3.4} by
$1-u$, we get
\begin{equation}\label{e-3.5}
   \frac{\Re[a_{n-1}]}{n}+\frac{\Re[a_{n-2}]}{n-1}(2+2u)\le\frac{1-\beta
   }{1-u}+\Ob{2}
\end{equation}
for each $\omega\ne1$.  In particular, inequality \eqref{e-3.5} holds 
for $u=u_1$ and $u=u_2$ as defined in Theorem \ref{1}.

Applying the linear transformation $\T$ defined in equation
\eqref{e-2.3} to inequality \eqref{e-3.5}, and using the values
computed in \eqref{e-2.4}, we see that
\begin{equation}\label{e-3.6}
   \Re[a_{n-1}]+\Re[a_{n-2}] \le(n+n D_1+D_2) (1-\beta)+\Ob{2}.
\end{equation}

Recall that 
$P'(z)=\prod_{j=1}^n (z-\zeta_j)=z^n+a_{n-1}z^{n-1}+\dots+a_0$, that
each $a_{n-k}=\Ob{}$ and that $\Re[a_{n-3}]=\Ob{2}$.  Then 
\begin{equation}\label{e-3.7}
\begin{split}
   |P|_\beta^{2n}&=(\min_j|\beta-\zeta_j|)^{2n} \le 
   \prod_{j=1}^n |\beta-\zeta_j|^2= |P'(\beta)|^2 \\
              &=P'(\beta)\overline P'(\beta)=\beta^{2n}+
   2\Re[a_{n-1}]\beta^{2n-1}+ 2\Re[a_{n-2}]\beta^{2n-2}+ \Ob{2} \\
              &=1-2n(1-\beta)+2\Re[a_{n-1}]+2\Re[a_{n-2}]+\Ob{2} \\
            &=\big[1-(1-\beta)+
                  (\Re[a_{n-1}]+\Re[a_{n-2}])/n\big]^{2n}+ \Ob{2}
\end{split}
\end{equation}
and so using inequalities \eqref{e-3.7} and then \eqref{e-3.6} we have
\begin{equation}\label{e-3.8}
\begin{split}
   |P|_\beta &\le 1-(1-\beta)+(\Re[a_{n-1}]+\Re[a_{n-2}])/n+\Ob{2} \\
             &\le 1+ (D_1+D_2/n) (1-\beta)+\Ob{2}.
\end{split}
\end{equation}
This completes our first step.

Our second step will be to verify the hypotheses of part 2 of Proposition
\ref{10}, by showing that
\begin{align*}
a_{n-1} &= n(1+D_1+D_2)(1-\beta)+\Ob{\alpha}, \\
a_{n-2} &= -(n-1)D_2(1-\beta)+\Ob{\alpha}, \text{ and} \\
a_{n-k} &= \Ob{\alpha} \text{ for $k\ge3$} .
\end{align*}

Combining inequalities \eqref{e-3.1} and \eqref{e-3.8}, we see that
\begin{equation}\label{e-3.9}
   |P|_\beta=1+(D_1+D_2/n) (1-\beta)+\Ob{2}. 
\end{equation}

Since equation \eqref{e-3.8} is thus an equality, then so are equations 
\eqref{e-3.7} and \eqref{e-3.6}, and thus equation \eqref{e-3.5} for 
$u=u_i$ and equations \eqref{e-3.4} and ~\eqref{e-3.3} for
$\Re[\omega]=u_i$.

Since equation \eqref{e-3.5} is an equality for $u=u_i$, we can solve
the resulting linear system in the variables $\Re[a_{n-1}]$ and
$\Re[a_{n-2}]$ and get
\begin{align*}
\Re[a_{n-1}]&= \frac{-n(u_1+u_2)}{(1-u_1)(1-u_2)}(1-\beta)+\Ob{2} \\
              &=  n(1+D_1+D_2)(1-\beta)+\Ob{2} \quad\text{and} \\
\Re[a_{n-2}]&= \frac{n-1}{2(1-u_1)(1-u_2)}(1-\beta)+\Ob{2} \\
             &= -(n-1)D_2(1-\beta)+\Ob{2}. 
\end{align*}
Note that from Proposition \ref{9} we have that $\Re[a_{n-k}]=\Ob{2}$
for $k\ge3$, so we now have the correct real parts for our second
step.  Thus we need only show that each $\Im[a_{n-k}]=\Ob{\alpha}$.

Recalling the definitions of $u_1$ and $u_2$ in Theorem \ref{1}, we can
choose $\omega_1$ and $\omega_2$ to be $(n+1)$st roots of $1$ so that 
$\Re[\omega_i]=u_i$. For $\omega=\omega_i$, expanding the 
products in equality \eqref{e-3.3} and cancelling those terms of
equality ~\eqref{e-3.4} gives us
\begin{equation}\label{e-3.10}
   \frac{\Im[a_{n-1}]}{n} \Im[\omega_i]+ 
   \frac{\Im[a_{n-2}]}{n-1} \Im[\omega_i^2]+ 
   \frac{\Im[a_{n-3}]}{n-2} \Im[\omega_i^3]=\Ob{2}.
\end{equation}

Consider the case $i=1$.  Since $|\omega_1|=1$ and since by part 1 of
Lemma \ref{8} we have $-1/2\le u_1<1$ then $\Im[\omega_1]\ne0$.  Now
by Proposition \ref{9}, $\Im[a_{n-k}]=\Ob{3/2}$ for $k\ge2$, so
equation \eqref{e-3.10} implies that $\Im[a_{n-1}]=\Ob{3/2}$. If $n=3$
or $n=5$ then by definition $\alpha=3/2$ so this completes our second
step for those two values of $n$.

Assume then without loss of generality that $n\ne3, 5$.  Again by 
part 1 of Lemma~\ref{8} we have $-1<u_2<u_1<1$ so $\Im[\omega_i]\ne0$.  Thus 
we may divide equation \eqref{e-3.10} by $\Im[\omega_i]$ to obtain
\begin{equation}\label{e-3.11}
   \frac{\Im[a_{n-1}]}{n}+
   \frac{\Im[a_{n-2}]}{n-1} (2 u_i)+ 
   \frac{\Im[a_{n-3}]}{n-2} (4 u_i^2-1)=\Ob{2}.
\end{equation}

Now subtracting equality \eqref{e-3.11} with $i=2$ from equality 
\eqref{e-3.11} with $i=1$ and dividing by $2(u_1-u_2)$ produces
\begin{equation}\label{e-3.12}
   \frac{\Im[a_{n-2}]}{n-1}+ \frac{\Im[a_{n-3}]}{n-2}2(u_1+ u_2) =
   \Ob{2}.
\end{equation}

Since equation \eqref{e-3.7} is an equality, we have each 
$|\beta-\zeta_j|=|P|_\beta+\Ob{2}$.  Recall that 
$\Im[a_{n-1}]=\Ob{3/2}$ and that 
$|P|_\beta=1+(D_1+D_2/n)(1-\beta)+\Ob{2}$.  Then by part
1 of Lemma \ref{13} we have $\Im[a_{n-2}]=(-3/2)\Im[a_{n-3}]+\Ob{5/2}$, so 
substituting into \eqref{e-3.12} we have 
\begin{equation*}
   \Im[a_{n-3}]\left[\frac{-3/2\  }{n-1} + \frac{2(u_1+u_2)}{n-2}\right]= 
   \Ob{2}.
\end{equation*}
Now by part 2 of Lemma \ref{8}, we have $u_1+u_2<0$ so the quantity in
brackets is non-zero. Then $\Im[a_{n-3}]=\Ob{2}$, and so solving back
in equations \eqref{e-3.12} and \eqref{e-3.11} we find that
$\Im[a_{n-k}]=\Ob{2}$ for all~$k\le 3$.  Note that by
Proposition~\ref{9}, we have $a_{n-k}=\Ob{2}$ for all $k \ge 4$, and
so $\Im[a_{n-k}]=\Ob{2}$ for all $k$.  Since $n\ne3, 5$ then by
definition $\alpha=2$ and so this finishes the proof of our second
step.

We will now finish the proof of Proposition \ref{5}. Consider
only those roots $z$ of~$P$ such that the nearest $\omega$ has
$\Re[\omega]=u_i$.   In our second step, we verified the hypotheses of 
part~2 of Proposition \ref{10}, so we have
\begin{equation*}
   |z|^{2n+2}=1-2(n+1)\Re[R]+(n+1)(\Gamma_1+\Gamma_2u_i)(1-\beta)^2+ 
   \Ob{\alpha+1}.
\end{equation*}

Since $|z|\le1$, this implies that
\begin{equation*}
   -\Re[R]\le -\frac{\Gamma_1+\Gamma_2u_i}{2}(1-\beta)^2+ \Ob{\alpha+1}
\end{equation*}
and so from the definition of $R$ in Proposition \ref{10} we have 
\begin{multline*}
-\Re\left[ a_{n-1}\frac{\overline\omega-1}{n}
     + a_{n-2}\frac{\overline\omega^2-1}{n-1} +\dots+a_0(\omega-1) 
\right] \\ 
    \le (1-\beta)-\frac{\Gamma_1+\Gamma_2u_i}{2}(1-\beta)^2+
\Ob{\alpha+1}.
\end{multline*}
Since $\Re(\omegabar)=u_i$, this inequality is also valid when
$\omega$ is replaced by $\omegabar$.  Note that by Proposition \ref{9} we 
have $\Re[a_{n-k}]=\Ob{3}$ for $k\ge5$, so averaging these two
inequalities gives us 
\begin{multline}\label{e-3.13}
\frac{\Re[a_{n-1}]}{n}\Re[1-\omega]+\dots+
\frac{\Re[a_{n-4}]}{n-3}\Re[1-\omega^4] \\
\le(1-\beta)-\frac{\Gamma_1+\Gamma_2u_i}{2}(1-\beta)^2+ \Ob{\alpha+1}.
\end{multline}

Note that since $|\omega|=1$, then $\Re[\omega^2]=2u_i^2-1$, 
$\Re[\omega^3]=4u_i^3-3u_i$ and $\Re[\omega^4]=8u_i^4-8u_i^2+1$.
Dividing inequality \eqref{e-3.13} by $1-u_i$, we get
\begin{multline*}
  \frac{\Re[a_{n-1}]}{n}
+\frac{\Re[a_{n-2}]}{n-1} (2+2u_i)
+\frac{\Re[a_{n-3}]}{n-2}(1+4u_i+4u_i^2)
+\frac{\Re[a_{n-4}]}{n-3}(8u_i^2+8u_i^3) \\
    \le \frac{1-\beta}{1-u_i}- \frac{(\Gamma_1+\Gamma_2u_i)(1-\beta)^2}{
2(1-u_i)}+\Ob{\alpha+1}. 
\end{multline*}
Applying to this the linear transformation $\T$ defined in 
\eqref{e-2.3} and using the values computed in \eqref{e-2.4}, we 
get an inequality of the form 
\begin{equation}\label{e-3.14}
\begin{split}
   \Re[a_{n-1}] &+ \Re[a_{n-2}]+c_3\Re[a_{n-3}]+c_4\Re[a_{n-4}] \\
                &\le (n+nD_1+D_2)(1-\beta) \\
              &\qquad-\big[(\Gamma_1/2)(n+nD_1+D_2)+ 
                (\Gamma_2/2)(nD_1+D_2)\big](1-\beta)^2 \\
              &\qquad+ \Ob{\alpha+1},
\end{split}
\end{equation}
where $c_3=\T(1+4u+4u^2)/(n-2)$ and $c_4=\T(8u^2+8u^3)/(n-3)$. 

Define
\begin{multline}\label{e-3.15}
Q = (-\Gamma_1/2)(n+nD_1+D_2)- (\Gamma_2/2)(nD_1+D_2) \\
      -(n-1)(n-2)(1-c_3)D_2(1+D_1+D_2/ n).
\end{multline}

Recall from our second step that for all $n$ we have that
$\Im[a_{n-1}]=\Ob{3/2}$, and that
$\Re[a_{n-2}]=-(n-\nobreak1)D_2(1-\nobreak\beta)+\Ob{2}$, and that
each $|\zeta_j-\beta|=1+(D_1+D_2/n)(1-\beta)+\Ob{2}$.  Then by part 2
of Lemma \ref{13}, we have
\begin{equation*}
   \Re[a_{n-3}]+2\Re[a_{n-4}]=-(n-1)(n-2)D_2(1+D_1+D_2/n)(1-\beta)^2+ 
   \Ob{3}.
\end{equation*}

Adding $1-c_3$ times this to inequality \eqref{e-3.14} gives us
\begin{multline}\label{e-3.16}
\Re[a_{n-1}]+\Re[a_{n-2}]+\Re[a_{n-3}]+(2-2c_3+c_4)\Re[a_{n-4}] \\
        \le (n+nD_1+D_2)(1-\beta)+Q(1-\beta)^2+\Ob{\alpha+1}. 
\end{multline}

Note that Lemma \ref{11} implies that $c_3<1/2$ for $n\ne3$, 4, and 6 and 
that $c_4\ge0$ for all $n$.  Using the definition of $\T$ in 
\eqref{e-2.3}, we calculate that for $n=4$ we have $c_3=3/2$ 
and $c_4=4$, and for $n=6$ we have $c_3=0.729$ and $c_4=0.972$. 
Thus for all $n\ge4$ we have $1-2c_3+c_4>0$. Note also that by our 
second step and Lemma~\ref{14} we have
$\Re[a_{n-4}]\ge\delta(1-\beta)^2+\Ob{\alpha+1}$.  Since $\delta=0$ 
except when $n=5$, and for $n=5$ we calculate $c_3=1/3$ and $c_4=2$, 
then 
\begin{equation*}
\begin{split}
-(1-2c_3+c_4)\Re[a_{n-4}] &\le 
-(1-2c_3+c_4)\delta(1-\beta)^2+\Ob{\alpha+1} \\
       &=(-7\delta/3)(1-\beta)^2+\Ob{\alpha+1}.
\end{split}
\end{equation*}
Adding this to equation \eqref{e-3.16} gives us
\begin{multline}\label{e-3.17}
\Re[a_{n-1}+a_{n-2}+a_{n-3}+a_{n-4}] \\
    \le (n+nD_1+D_2)(1-\beta)+
(Q-7\delta/3)(1-\beta)^2+\Ob{\alpha+1}.
\end{multline}

Let 
\begin{align*}
Q_1&=-n(1-\beta)+a_{n-1}+a_{n-2}+a_{n-3}+a_{n-4}+a_{n-5} 
\quad\text{and} \\
Q_2&=n(n-1)(1-\beta)^2/2-\big[(n-1)a_{n-1}+(n-2)a_{n-2}\big](1-\beta).
\end{align*}

Recall from our first step that each $a_{n-k}=\Ob{}$ so $Q_1=\Ob{}$ 
and $Q_2=\Ob{2}$.

Now from our second step we know that $a_{n-k}=\Ob{\alpha}$ for 
$k\ge3$, and from Proposition \ref{9} we know that $a_{n-k}=\Ob{3}$
for $k\ge6$, so 
\begin{equation*}
\begin{split}
P'(\beta)&=\beta^n+a_{n-1}\beta^{n-1}+\dots+a_0 \\
         &=1-n(1-\beta)+\frac{n(n-1)}{2}(1-\beta)^2 
+a_{n-1}\big[1-(n-1)(1-\beta)\big] \\
         &\qquad +a_{n-2}\big[1-(n-2)(1-\beta)\big]
+a_{n-3}+a_{n-4}+a_{n-5}+\Ob{\alpha+1} \\
         &=1+Q_1+Q_2+\Ob{\alpha+1}.
\end{split}
\end{equation*}

Then $|P'(\beta)|^2=P'(\beta)\overline{P'(\beta)}=
1+2\Re[Q_1]+2\Re[Q_2] +|Q_1|^2+\Ob{\alpha+1}$.  Note from our second 
step that each $\Im[a_{n-k}]=\Ob{\alpha}$ so $\Im[Q_1]=\Ob{\alpha}$. 
Then $(1+\Re[Q_1]+\Re[Q_2])^2= |P'(\beta)|^2+\Ob{\alpha+1}$ and so 
$|P'(\beta)|=1+\Re[Q_1]+\Re[Q_2]+\Ob{\alpha+1}$.  Substituting the 
values of $Q_1$ and $Q_2$ and using the results of our second step 
gives us 
\begin{equation}\label{e-3.18}
\begin{split}
   |P'(\beta)| &=1-n(1-\beta)+\Re[a_{n-1}+a_{n-2}+a_{n-3}+a_{n-4}] \\
           &\qquad +(n-1)\big[n/2-n(1+D_1+D_2)+
   (n-2)D_2\big](1-\beta)^2 \\
          &\qquad+\Ob{\alpha+1}.
\end{split} 
\end{equation}

Using the first line of inequality \eqref{e-3.7}, then inequalities 
\eqref{e-3.18} and \eqref{e-3.17}, we have
\begin{equation}\label{e-3.19}
\begin{split}
   |P|_\beta^n &\le |P'(\beta)| \\ 
         &\le 1+(nD_1+D_2)(1-\beta) \\
         &\qquad +\big[Q-7\delta/3
   -(n-1)(n/2+nD_1+2D_2)\big](1-\beta)^2 \\
         &\qquad+ \Ob{\alpha+1}.
\end{split}
\end{equation}

We seek now to compute the coefficient of $(1-\beta)^2$ in this 
inequality.  Note first that from the definitions of $\Gamma_1$ and 
$\Gamma_2$ in Proposition \ref{10} we have 
\begin{equation*}
\begin{split}
-\frac{\Gamma_1}{2}(n+nD_1+D_2)&-\frac{\Gamma_2}{2}(nD_1+D_2) \\
  & = -\frac{\Gamma_1+\Gamma_2}{2}(n+nD_1+D_2)+\frac{n\Gamma_2}{2} \\
  & =\left[ (1+2D_1)n-\left(\frac{1}{2}+2D_1-2D_2\right) \right] 
\big[(1+D_1)n+D_2
\big] \\
   &\qquad+ n(1+D_1+D_2) \big[nD_2+(D_1-2D_2)\big].
\end{split}
\end{equation*}

Now from the definition of $c_3$ (after inequality \eqref{e-3.14}) 
combined with equalities~\eqref{e-2.4} we have
$(n-2)c_3D_2=-(n+1+D_1+3nD_1+3D_2)$ and so
\begin{equation*}
   (n-2)(1-c_3)D_2=(1+3D_1+D_2)n+(1+D_1+D_2).
\end{equation*}
Substituting these values into equation \eqref{e-3.15} and collecting 
like powers of $n$, we conclude that
\begin{multline}
Q=\big[-D_1-D_1^2+D_2^2\big]n^2+ 
          \left[-\frac{1}{2}+\frac{1}{2}D_1+D_1^2-3D_2^2\right]n \\
  + \left[1+2D_1+\frac{1}{2}D_2+D_1^2+D_1D_2+2D_2^2\right] 
       +\big[D_2+D_1D_2+D_2^2\big]/n
\end{multline}
and so comparing this with the definition of $D$ in Theorem \ref{1}, we see 
that 
\begin{equation}\label{e-3.20}
   Q-(n-1)(n/2+nD_1+2D_2)=nD+\frac{n(n-1)}{2}(D_1+D_2/n)^2.
\end{equation}
Substituting this into inequality \eqref{e-3.19}, we have
\begin{equation*}
\begin{split}
|P|_\beta^n &\le 1+(nD_1+D_2)(1-\beta) \\
   &\qquad+ \left[nD+\frac{n(n-1)}{2}(D_1+D_2/n)^2
-7\delta/3\right](1-\beta)^2
         +\Ob{\alpha+1} \\
   &=\left[1+(D_1+D_2/n)(1-\beta)+ 
       \left(D-\frac{7\delta}{3n}\right)(1-\beta)^2\right]^n+ 
\Ob{\alpha+1}.
\end{split} 
\end{equation*}
Note that (from the definitions of $\delta$ in Lemma \ref{14} and $\Delta$ in
Theorem \ref{1}) for all $n$ we have $\Delta=-7\delta/(3n)$, and so 
\begin{equation*}
   |P|_\beta \le 1+(D_1+D_2/n)(1-\beta)+ 
   (D+\Delta)(1-\beta)^2+\Ob{\alpha+1}. 
\end{equation*}
This completes the proof of Proposition \ref{5}.


\section{Proof of Proposition \ref{6}}

This proof parallels the proof of \cite[Theorem 2]{Mil-3}.   We begin 
by letting 
\begin{align*}
u &= \frac{-i\sqrt{15}}{15}(1-\beta)^{1/2}-\frac{6}{10}(1-\beta)+
\frac{i\sqrt{15}}{300}(1-\beta)^{3/2}- \frac{33}{600}(1-\beta)^2 \\
\noalign{\text{\noindent and}}
v &= 
\frac{4i\sqrt{15}}{15}(1-\beta)^{1/2}-\frac{1}{10}(1-\beta)+ 
\frac{46i\sqrt{15}}{300}(1-\beta)^{3/2}+ \frac{532}{600}(1-\beta)^2.
\end{align*}
Let $P'(z)=(z-u)^4(z-v)$ and let $P(z)=\int_\beta^zP'(t)\,dt$.  Note 
that $u-\beta=-1+u+(1-\beta)$ so
\begin{equation*}
\begin{split}
|u-\beta|^2 &= \big[-1+(4/10)(1-\beta)-(33/600)(1-\beta)^2\big]^2 \\
     &\qquad+ \big[(-\sqrt{15}/15)(1-\beta)^{1/2}+
(\sqrt{15}/300)(1-\beta)^{3/2}\big]^2 \\
               &= 1-(11/15)(1-\beta)+ (79/300)(1-\beta)^2+\Ob{3}
\end{split}
\end{equation*}
and $v-\beta=-1+v+(1-\beta)$ so
\begin{equation*}
\begin{split}
|v-\beta|^2&=\big[-1+(9/10)(1-\beta)+(532/600)(1-\beta)^2\big]^2 \\
    &\qquad + \big[(4\sqrt{15}/15)(1-\beta)^{1/2}+
(46\sqrt{15}/300)(1-\beta)^{3/2}\big]^2 \\
          &=1-(11/15)(1-\beta)+(79/300)(1-\beta)^2+\Ob{3}.
\end{split}
\end{equation*}
Now 
\begin{multline*}
\big[1-(11/30)(1-\beta)+(29/450)(1-\beta)^2\big]^2 \\
=1-(11/15)(1-\beta)+(79/300)(1-\beta)^2+\Ob{3},
\end{multline*}
and so we have
\begin{equation*}
\begin{split}
|P|_\beta &= \min\{|u-\beta|,|v-\beta|\} \\
          &= 1-(11/30)(1-\beta)+ (29/450)(1-\beta)^2+\Ob{3}.
\end{split}
\end{equation*}

By definition $P$ is of degree 6 and $P(\beta)=0$.  Thus to verify 
that $P \in S(6,\beta)$ we need only show that all the roots of $P$
remain in the closed unit disk when $\beta$ is sufficiently close to
1.  Now
\begin{align*}
u^2 &= (-1/15)(1-\beta)+ (2i\sqrt{15}/25)(1-\beta)^{3/2}+\Ob{2}, \\
u^3 &= (i\sqrt{15}/225)(1-\beta)^{3/2}+\Ob{2},\quad\text{and} \\
u^4 &= \Ob{2},
\end{align*}
so writing $P'(z)=z^5+a_4z^4+\dots+a_0$, we calculate that
\begin{align*}
a_4 &= -(4u+v)=(5/2)(1-\beta)-(i\sqrt{15}/6)(1-\beta)^{3/2}- 
           (2/3)(1-\beta)^2 \\ 
a_3 &= u(6u+4v) \\
     &= (2/3)(1-\beta)- (2i\sqrt{15}/15)(1-\beta)^{3/2}+ 
          3(1-\beta)^2+ \Ob{5/2} \\
a_2 &= -u^2(4u+6v)=(4i\sqrt{15}/45)(1-\beta)^{3/2}+
          (7/5)(1-\beta)^2+\Ob{5/2} \\
a_1 &= u^3(u+4v) = (-1/15)(1-\beta)^2+\Ob{5/2} \\
a_0 &= -u^4v = \Ob{5/2}.
\end{align*}

Recall from the values computed at the beginning of section 2 that for
$n=5$ we have $\alpha=3/2$, $u_1=-1/2$, $u_2=-1$, $D_1=-1/3$ and
$D_2=-1/6$.  Note that in part 2 of Proposition \ref{10} the values of
the $a_k$'s computed above satisfy the hypotheses, and that
$\Gamma_2=-5/6$ and $\Gamma_1=-13/6$.

Let us apply part 2 of Proposition \ref{10}  to the case $\omega=-1$.  Note
that $\Re[\omega]=u_2$ and $\Gamma_1+\Gamma_2u_2=-4/3$.  Since 
$\omega=-1$ we have
\begin{equation*}
   R = (1-\beta)-(2/5)a_4-(2/3)a_2-2a_0,
\end{equation*}
and so 
\begin{equation*}
\begin{split}
\Re[R] &= (1-\beta)- (2/5)\big[(5/2)(1-\beta)- 
(2/3)(1-\beta)^2\big] \\
    &\qquad- (2/3)(7/5)(1-\beta)^2+\Ob{5/2} \\
     &=(-2/3)(1-\beta)^2+\Ob{5/2}.
\end{split}
\end{equation*}
Thus by part 2 of Proposition \ref{10} we have 
\begin{equation*}
\begin{split}
|z|^{12} &= 1-12(-2/3)(1-\beta)^2+6(-4/3)(1-\beta)^2+\Ob{5/2} \\
     &= 1+\Ob{5/2},
\end{split}
\end{equation*}
and so $|z|=1+\Ob{5/2}$.

Let us now apply part 2 of Proposition \ref{10} to the case
$\omega=(1/2)(-1\pm i\sqrt3)$.  Note that $\Re[\omega]=u_1$ and 
$\Gamma_1+\Gamma_2u_1=-7/4$.  Now
\begin{multline*}
R=(1-\beta)+(a_4/10)(-3\mp i\sqrt3)+ (a_3/8)(-3\pm i\sqrt3) \\
         +(a_1/4)(-3\mp i\sqrt3)+ (a_0/2)(-3\pm i\sqrt3)
\end{multline*}
so
\begin{equation*}
\begin{split}
\Re[R] &= (1-\beta)- 
(3/10)\big[(5/2)(1-\beta)-(2/3)(1-\beta)^2\big] \\
     &\qquad\pm (\sqrt3/10)(-\sqrt{15}/6)(1-\beta)^{3/2}- 
(3/8)\big[(2/3)(1-\beta)+3(1-\beta)^2\big] \\
     &\qquad \mp (\sqrt3/8)(-2\sqrt{15}/15)(1-\beta)^{3/2}- 
(3/4)(-1/15)(1-\beta)^2+\Ob{5/2} \\
     &=(-7/8)(1-\beta)^2+\Ob{5/2}.
\end{split}
\end{equation*}
Thus by part 2 of Proposition \ref{10} we have 
\begin{equation*}
\begin{split}
|z|^{12} &= 1-12(-7/8)(1-\beta)^2+6(-7/4)(1-\beta)^2+\Ob{5/2} \\
          &=1+\Ob{5/2},
\end{split}
\end{equation*}
so $|z|=1+\Ob{5/2}$.

Finally, let us apply part 1 of Proposition \ref{10} with $r=1$ to the case 
$\omega=(1/2)(1\pm\nobreak i\sqrt3)$.   Note that 
\begin{equation*}
   R = (1-\beta)+(a_4/10)(-1\mp i\sqrt3)+ (a_3/8)(-3\mp 
   i\sqrt3)+\Ob{3/2}
\end{equation*}
so 
\begin{equation*}
\begin{split}
\Re[R]&=(1-\beta)+(-1/10)(5/2)(1-\beta)+
(-3/8)(2/3)(1-\beta)+\Ob{3/2} \\
                 &=(1/2)(1-\beta)+\Ob{3/2}.
\end{split}
\end{equation*}
Thus by part 1 of Proposition \ref{10} we have $|z|^2=1-(1-\beta)+\Ob{3/2}$ 
and so $|z|=1-(1/2)(1-\beta)+\Ob{3/2}$.

At this stage, we know that $|P|_\beta= 1-(11/30)(1-\beta)+
(29/450)(1-\beta)^2+\Ob{3}$ and that if $\beta$ is sufficiently close
to $1$ then all roots $z$ of $P$ have $|z|\le1+\Ob{5/2}$.  Since the
roots of $P$ approach the roots of $z^6-1$, then the non-$\beta$ roots
of $P$ are bounded away from $\beta$.  Thus by Lemma \ref{12}, there
is a polynomial $Q \in S(6,\beta)$ with $|Q|_\beta=
1-(11/30)(1-\beta)+ (29/450)(1-\beta)^2+\Ob{5/2}$.  This completes the
proof of Proposition \ref{6}.


\section{Proof of Proposition \ref{7}}

Let $b_1=1+D_1+D_2/n$, let $b_2=(n-1)D_2$, and let
$z_0=-b_1(1-\beta)-D(1-\beta)^2$. Then 
$z_0-\beta=-1+(1-b_1)(1-\beta) -D(1-\beta)^2$, and (for $\beta$
near $1$) this is real and negative so $|z_0-\beta|=1+
(D_1+D_2/n)(1-\beta)+D(1-\beta)^2$.

Now let $x$ be a real constant, depending only on $n$ (and to be 
determined later), and let 
\begin{multline*}
q(z)=z^2+ \big[(b_2+2b_1)(1-\beta)-2x(1-\beta)^2 \big] z \\
  +\big[-b_2(1-\beta)+ (b_1^2+b_2+2D+2x)(1-\beta)^2\big].
\end{multline*}
Now by part 4 of Lemma \ref{8} we have $D_2<0$ and so $b_2<0$.  Since the 
discriminant of $q(z)$ is $4b_2(1-\beta)+\Ob{2}$, then (for $\beta$ 
near $1$) the roots of $q$ are complex conjugates.  If we denote these
roots by $z_1$ and $\overline z_1$ then by writing $\beta=1-(1-\beta)$
we have
\begin{equation*}
\begin{split}
|z_1-\beta|^2 &= (z_1-\beta)(\overline z_1-\beta)=q(\beta) \\
            &= 1+(2b_1-2)(1-\beta)+(1-2b_1+b_1^2+2D)(1-\beta)^2 
+\Ob{3} \\ 
        &=\big[1+(b_1-1)(1-\beta)+D(1-\beta)^2\big]^2+\Ob{3}
\end{split}
\end{equation*}
so $|z_1-\beta|=1+(D_1+D_2/n)(1-\beta)+D(1-\beta)^2+\Ob{3}$.

Let $P'(z)=(z-z_0)^{n-2}q(z)$ and $P(z)=\int_\beta^z P'(t)\,dt$, so 
\begin{equation*}
   |P|_\beta=1+(D_1+D_2/n)(1-\beta)+D(1-\beta)^2+\Ob{3}.
\end{equation*}

Now $z_0=\Ob{}$ so
\begin{equation*}
\begin{split}
(z-z_0)^{n-2} &= z^{n-2}-(n-2)z_0z^{n-3}+\binom{n-2}2 
z_0^2z^{n-4}+\Ob{3} \\
          &= z^{n-2}+(n-2)\big[b_1(1-\beta) 
                +D(1-\beta)^2\big]z^{n-3} \\
          &\qquad+ \binom{n-2}2 b_1^2(1-\beta)^2z^{n-4}+\Ob{3}.
\end{split}
\end{equation*}
Then letting $t_1=(n^2-n)b_1^2/2+(n-2)b_1b_2+b_2$ we have
\begin{equation}\label{e-5.1}
\begin{split} 
   P'(z) &= (z-z_0)^{n-2}q(z) \\
        &= z^n+ \big[(nb_1+b_2)(1-\beta)+ 
                (nD-2D-2x)(1-\beta)^2\big]z^{n-1} \\
        &\qquad+ \big[-b_2(1-\beta)+ 
                  (t_1+2D+2x)(1-\beta)^2\big]z^{n-2} \\
        &\qquad-(n-2)b_1b_2(1-\beta)^2z^{n-3}+\Ob{3}.
\end{split}
\end{equation}

Note that by its definition, $P$ is a polynomial of degree $n+1$ and 
$P(\beta)=0$.  Thus to show that $P\in S(n+1,\beta)$ it will suffice 
to show that all roots of $P$ remain in the unit disk when $\beta$ is
sufficiently close to $1$.

Let $\omega\ne1$ be an $(n+1)$th root of $1$, let $u=\Re[\omega]$ and 
note that since $|\omega|=1$ then $\Re[\omega^2]=2u^2-1$,
$\Re[\omega^3]=4u^3-3u$, and $\omega^{n-k}=\omegabar^{k+1}$. 
Substituting the coefficients of equation \eqref{e-5.1} into the
formula for $R$ in Proposition \ref{10}, we have 
\begin{multline*}
R=(1-\beta)+(nb_1+b_2)(1-\beta)(\omegabar-1)/n \\
     -b_2(1-\beta)(\omegabar^2-1)/(n-1)+\Ob{2}. 
\end{multline*}
Substituting the values of $b_1$ and $b_2$ into this formula, we see 
by part 1 of Proposition~\ref{10} with $r=1$ that 
\begin{equation*}
   |z|^2=1-2(1-\beta)\big[1+(1+D_1+D_2)(u-1)-D_2(2u^2-2)\big] +\Ob{2}.
\end{equation*}

Recall from part 4 of Lemma \ref{8} that $D_2<0$, so the quantity in
square brackets is quadratic in $u$ with positive leading coefficient.
By elementary calculus, its minimum (over all real numbers) occurs
when $1+D_1+D_2-4D_2u=0$, which happens when $u=(1+D_1+D_2)/(4D_2)=
(u_1+u_2)/2$, which is between $u_1$ and $u_2$. Now $u_1$ and $u_2$
are (by definition) the real parts of adjacent $(n+1)$th roots of $1$,
so there are no possible values of $u$ between $u_1$ and $u_2$, so the
minimum (over all possible values of $u$) must occur at either $u_1$
or $u_2$.  From part 7 of Lemma \ref{8} we see that at these values
the quantity in square brackets is $0$, and so the minimum value of
the quantity in square brackets is 0.  Thus for $\Re[\omega]\ne u_i$
the quantity in square brackets is positive, so for these values of
$\omega$ and for $\beta$ sufficiently close to ~$1$ we have $|z|<1$,
and so these roots remain in the unit disk.

Thus we need only concern ourselves with the case $\Re[\omega]=u_i$. 
In this case,  by part 2 of Proposition \ref{10} we have
\begin{equation*}
   |z|^{2n+2}=1-2(n+1)\Re[R]+(n+1)(\Gamma_1+\Gamma_2u_i)(1-\beta)^2+ 
   \Ob{\alpha+1}. 
\end{equation*}
To get $P \in S(n+1,\beta)$ we will seek a value of $x$ so that
$|z|=1+\Ob{\alpha+1}$, so we will need
\begin{equation}\label{e-5.2}
   \Re[R]-(1/2)(\Gamma_1+\Gamma_2u_i)(1-\beta)^2=\Ob{\alpha+1}
\end{equation}
for both $i=1$ and $i=2$.  

Substituting the coefficients of equation \eqref{e-5.1} into the
formula for $R$ in Proposition \ref{10}, we have
\begin{equation}\label{e-5.3}
\begin{split}
   R &= (1-\beta)+\big[(nb_1+b_2)(1-\beta)+ (nD-2D-2x)(1-\beta)^2\big] 
   (\omegabar-1)/n \\
     &\qquad+\big[-b_2(1-\beta)+ (t_1+2D+2x)(1-\beta)^2 \big] 
   (\omegabar^2-1)/(n-1) \\
    &\qquad-(n-2)b_1b_2(1-\beta)^2(\omegabar^3-1)/(n-2)+\Ob{3}.
\end{split}
\end{equation}
Taking the real parts of equation \eqref{e-5.3} and collecting like 
powers of $(1-\beta)$ gives us
\begin{equation*}
\begin{split}
\Re[R] &= \big[1+(nb_1+b_2)(u_i-1)/n-
          b_2(2u_i^2-2)/(n-1)\big](1-\beta) \\
       &\qquad+ \bigg[(nD-2D-2x)(u_i-1)/n+ 
(t_1+2D+2x)(2u_i^2-2)/(n-1) \\
       &\qquad\qquad -b_1b_2(4u_i^3-3u_i-1)\bigg](1-\beta)^2 +\Ob{3}.
\end{split}
\end{equation*}

Substituting the values of $b_1$ and $b_2$ into this formula, we see 
from part 7 of Lemma \ref{8} that the coefficient of $(1-\beta)$ in $\Re[R]$
is zero, so to satisfy equation \eqref{e-5.2} we need only find a value
of $x$ such that the coefficient of $(1-\beta)^2$ in equation
\eqref{e-5.2} is 0. We divide this coefficient by $u_i-1$ and denote 
the result by $Z_i$, so
\begin{multline}
Z_i=(nD-2D-2x)/n+(t_1+2D+2x)(2u_i+2)/(n-1) \\
  -(n-1)D_2(1+D_1+D_2/n)(4u_ i^2+4u_i+1)+
(1/2)(\Gamma_1+\Gamma_2u_i)/(1-u_i).
\end{multline}

Note that the coefficient of $x$ in $Z_i$ is $-2/n+(4u_i+4)/(n-1)$, 
which is nonzero by part 3 of Lemma \ref{8}, so each equation $Z_i=0$ has a 
solution for $x$.  To show that these solutions are identical, we 
will show that $Z_1$ and $Z_2$ (considered as linear expressions in 
the variable $x$) are scalar multiples of each other.

To see this, we eliminate $x$ by applying the transformation $\T$ 
defined in equation~\eqref{e-2.3}. Since in equation \eqref{e-3.14} we 
defined $c_3=\T(1+4u+4u^2)/(n-2)$, then from equations \eqref{e-2.4} we
see that 
\begin{multline}
\T(Z_i)=nD+t_1-(n-1)(n-2)c_3D_2(1+D_1+D_2/n) \\
+(\Gamma_1/2)(n+nD_1 + D _2)+(\Gamma_2/2)(nD_1+D_2).
\end{multline}
Comparing this to the value of $Q$ defined in equation \eqref{e-3.15}, 
we see that
\begin{equation}\label{e-5.4}
   \T(Z_i)=nD+t_1-Q-(n-1)(n-2)D_2(1+D_1+D_2/n).
\end{equation}

Note that by equation \eqref{e-3.20} we have
\begin{equation*}
   Q=nD+\frac{n(n-1)}{2}(D_1+D_2/n)^2+(n-1)(n/2+nD_1+2D_2).
\end{equation*}
Substituting the values of $b_1$ and $b_2$ into our definition of 
$t_1$ gives us
\begin{equation*}
   t_1=(n-1)\big[(n/2)(1+D_1+D_2/n)^2+(n-2)D_2(1+D_1+D_2/n)+D_2\big]
\end{equation*}
and so $Q-t_1=nD-(n-1)(n-2)D_2(1+D_1+D_2/n)$.  Substituting this into 
equation \eqref{e-5.4} gives us $\T(Z_i)=0$.  Since $\T(Z_i)$ is a 
linear combination of $Z_1$ and $Z_2$, this implies that $Z_1$ and 
$Z_2$ (considered as polynomials in $x$) are scalar multiples of one 
another, and so there is a single value of $x$ that satisfies 
equation \eqref{e-5.2} for both $i=1$ and $i=2$.

Using this value of $x$, we have now constructed a real polynomial $P$ 
with 
\begin{equation*}
   |P|_\beta=1+(D_1+D_2/n)(1-\beta)+ D(1-\beta)^2+\Ob{3}
\end{equation*}
and such that all roots $z$ of $P$ have $|z|\le1+\Ob{\alpha+1}$.
Since the roots of $P$ approach the roots of $z^{n+1}-1$, then the
non-$\beta$ roots of $P$ are bounded away from~$\beta$.  Thus by Lemma
\ref{12}, there is a real polynomial $Q \in S(n+1,\beta)$ with
\begin{equation*}
   |Q|_\beta=1+(D_1+D_2/n)(1-\beta)+ D(1-\beta)^2+\Ob{\alpha+1}.
\end{equation*}
This finishes the proof of Proposition \ref{7}.


\end{document}